\documentclass[10pt,a4paper]{article}
\usepackage{amsmath}
\usepackage{amsfonts}
\usepackage{amsthm}
\usepackage{amssymb}
\usepackage{enumerate}
\theoremstyle{plain}
\newtheorem{theo}{Theorem}[section]%[Definizioni]
\newtheorem{lem}[theo]{Lemma}%[Definizioni]
\newtheorem{cor}[theo]{Corollary}%

\theoremstyle{definition}
\newtheorem{definition}[theo]{Definition}%[Definizioni]

\theoremstyle{remark}
\newtheorem{ex}[theo]{Example}%[Definizioni]
\newtheorem{rem}[theo]{Remark}%[Definizioni]

\numberwithin{equation}{section}

\newcommand{\R}{\mathbb{R}}

\newcommand{\M}{\mathbb{M}}

\title{The linear constraints \\ in Poincar\'{e} and Korn type inequalities
\thanks{Work supported in part by MIUR, PRIN n. 2004011204.}}
\author{Giovanni Alessandrini\thanks{Dipartimento di Matematica e Informatica,
Universit\`a degli Studi di Trieste, Italy,
\textsf{alessang@univ.trieste.it}}\ , Antonino
Morassi\thanks{Dipartimento di Georisorse e Territorio,
Universit\`a degli Studi di Udine, Italy,
\textsf{antonino.morassi@uniud.it}}\ \ and Edi
Rosset\thanks{Dipartimento di Matematica e Informatica,
Universit\`a degli Studi di Trieste, Italy,
\textsf{rossedi@univ.trieste.it}}}
\date{}

\begin{document}

\maketitle

\begin{abstract}
We investigate the character of the linear constraints which are
needed for Poincar\'e and Korn type inequalities to hold. We
especially analyze constraints which depend on restriction on
subsets of positive measure and on the trace on a portion of the
boundary.
\end{abstract}

\section{Introduction} \label{sec:
introduction} Let us consider a bounded domain $\Omega\subset
\R^n$ with Lipschitz boundary. The nature of Poincar\'e inequality
is that a function $u$ (say in $W^{1,2}(\Omega)$) is uniquely
determined by its gradient $\nabla u\in L^2(\Omega)$ up to an
additive constant.

Such a constant is usually chosen to be the average of $u$ in
$\Omega$, but also the average of its trace on $\partial\Omega$
would do for this purpose. Similarly if one takes a weighted
average on suitable subsets  $E$ of $\overline{\Omega}$. Meyers
\cite{l:mey} has examined this subject in great depth
characterizing the sets $E \subset \overline{\Omega}$ for which
this is possible in terms of their capacity. His analysis in fact
extends to Poincar\'e inequalities in $W^{m,p}(\Omega)$, for any
$m=1,2, \ldots$ and $1<p<\infty$. See also Meyers and Ziemer
\cite{l:meyz} for the extreme case of $BV(\Omega)$ and Ziemer
\cite[Chapter 4]{l:z} for a more recent account on such results.

In the case of Korn inequality, a vector valued function $u\in
W^{1,2}(\Omega;\R^n)$ is uniquely determined by the strain tensor
$\frac{1}{2}(\nabla u +(\nabla u)^T)$ up to an infinitesimal rigid
displacement $r(x)=b+Mx$, that is an affine transformation for
which the matrix $M=\nabla r$ is skew symmetric. In this case,
provided the origin of the coordinates is placed in the center of
mass of $\Omega$, $x_\Omega=\frac{1}{|\Omega|}\int_\Omega x$, an
obvious choice is given by $b=\frac{1}{|\Omega|}\int_\Omega u$,
$M=\frac{1}{|\Omega|}\int_\Omega \frac{1}{2}(\nabla u -(\nabla
u)^T)$.

In other words, the common feature of Poincar\'e and Korn
inequalities is that they  are inequalities in which a norm is
dominated by a seminorm for functions which satisfy a suitable
finite set of linear constraints.

In this paper we examine the character of such constraints and in
particular we investigate whether such constraints can be chosen
in such a way that they can depend only on the values of the
function $u$ restricted to a subset of $\Omega$ or of its boundary
$\partial\Omega$. We are especially interested in obtaining
constructive evaluations of the constants in a number of specific
cases. In fact, constraints of such type are useful in many
instances and in particular, in recent years, they have shown up
in connection with inverse boundary value problems, for instance
in estimates of unique continuation for the system of elasticity,
\cite{l:amr2}, \cite{l:mr}. For such purposes, it is mandatory to
have a concrete and constructive evaluation of the constants
involved. In this note we develop this aspect, improving some
known estimates, and obtaining some which were not available in
the literature.

We deal with three types of inequalities, Poincar\'e inequality in
$W^{1,2}(\Omega)$, Poincar\'e inequality in $W^{2,2}(\Omega)$ and
Korn inequality in $W^{1,2}(\Omega)$. The results are all based on
an elementary functional analytic argument (Lemma \ref{lem:2.1})
which is due to N. Meyers \cite[Proposition 1]{l:mey}. And indeed
the consequences  for the standard Poincar\'e inequality, which we
summarize in Section 3 are also found in \cite{l:mey}.

This is not completely the case for the Poincar\'e inequality in
$W^{2,2}(\Omega)$, Section 4, and Korn inequality, Section 5.

In Section 4 we treat Poincar\'e inequality in $W^{2,2}(\Omega)$,
deriving, along the lines of Meyers' result, the concrete
evaluations of the constants in various instances, see Examples
\ref{ex:Poincare_E_H2}, \ref{ex:Poincare_balls_H2}. We observe
however that, when one deals with an open portion $\Gamma$ of the
boundary $\partial\Omega$, the suitable linear constraints that
occur in Meyers' approach, require the knowledge of the boundary
trace of $u$ \emph{and} of its first derivatives. We thus examine
the possibility to determine linear constraints, suitable for the
validity of a Poincar\'e type inequality, which \emph{only} depend
on the trace of $u$ on $\Gamma$. See Example
\ref{ex:Poincare_Gamma_H2}, and in particular Theorem
\ref{theo:funzionale_traccia}, for the Poincar\'e inequality in
$W^{2,2}(\Omega)$. It is peculiar in this case that one has to
distinguish the cases when $\Gamma$ is a subset of an
$(n-1)$-dimensional hyperplane or not.

In Section 5 we treat Korn inequality and obtain some estimates
that we believe to be new.
See Examples \ref{ex:korn_E} and \ref{ex:korn_balls}.
In particular, inequality \eqref{eq:5.korn_E} generalizes
a classical inequality proven by
Kondrat'ev and Oleinik (\cite[Theorem 1]{l:ko}) when
$E$ is a ball and $\Omega$ is starlike with respect to $E$.
We also examine the case of
constraints which depend on the trace of $u$ on an open
portion $\Gamma$ of $\partial\Omega$. In this case, we find
that the constraints can be
chosen to depend on the trace of $u$ on $\Gamma$ for any
open subset $\Gamma$ of $\partial\Omega$, see
Theorem \ref{theo:funzionale_traccia_korn}
and Corollary
\ref{cor:corollario_funzionale_traccia_korn}, which in fact
substantially generalize a
classical inequality proven by Kondratev and Oleinik
(\cite[Theorem 2]{l:ko}) for functions $u$ vanishing on a basis
of a cylinder $\Omega$.

\section{A basic Lemma} \label{sec:
Lemma}

\begin{lem}
  \label{lem:2.1}
Let $X$, $Y$ be Banach spaces, and let $L:X\rightarrow Y$ a
bounded linear operator. Let $X_0\subset X$ be its null-space and
let $P:X\rightarrow X_0$ be a bounded linear operator such that
$P_{|X_0}=Id$. Assume that there exists $K>0$ such that
\begin{equation}
  \label{eq:2.1}
  \|x-Px\|_X\leq K\|Lx\|_Y \ , \quad\hbox{for every } x\in X \ .
\end{equation}
For every bounded linear operator $T:X\rightarrow X_0\subset X$
such that $T_{|X_0}=Id$ we have
\begin{equation}
  \label{eq:2.2}
  \|x-Tx\|_X\leq (1+\|T\|)K\|Lx\|_Y \ , \quad\hbox{for every } x\in X \ .
\end{equation}
Here $\|T\|$ denotes the operator norm on ${\cal L}(X,X)$.
\end{lem}
\begin{proof}
The argument is due to N. Meyers \cite[Proposition 1]{l:mey}, we
reproduce it here, because of its brevity. We have
\begin{equation*}
 x-Tx=x-Px-Tx+Px=x-Px-T(x-Px) \ .
\end{equation*}
Hence
\begin{equation*}
  \|x-Tx\|_X\leq (1+\|T\|)\|x-Px\|_X\leq (1+\|T\|)K\|Lx\|_Y \ .
\end{equation*}
\end{proof}
\begin{rem}
  \label{rem:2.2}
Note that the condition $T_{|X_0}=Id$ is in fact necessary if an
inequality of the following form holds true
\begin{equation*}
  \|x-Tx\|_X\leq Const.\|Lx\|_Y \ .
\end{equation*}
Indeed, if $x\in X_0$, then $Lx=0$ and hence $Tx=x$. Note also
that the evaluation of the constant in \eqref{eq:2.2} is not
optimal, as is evident if one chooses $T=P$. When $T$ is close to
$P$ the following evaluation may be more convenient. Since
$(T-P)Px=0$, we have
\begin{equation*}
 x-Tx=x-Px-(T-P)x=x-Px-(T-P)(x-Px) \ ,
\end{equation*}
hence
\begin{equation}
  \label{eq:2.3}
  \|x-Tx\|_X\leq (1+\|T-P\|)K\|Lx\|_Y \ , \quad\hbox{for every } x\in X \ .
\end{equation}
Therefore \eqref{eq:2.2} can be improved as follows
\begin{equation}
  \label{eq:2.4}
  \|x-Tx\|_X\leq (1+\min\{\|T\|,\|T-P\|\})K\|Lx\|_Y \ , \quad\hbox{for every } x\in
X \ .
\end{equation}
\end{rem}

\section{The Poincar\'e inequality}
\label{sec: Poincare}

We shall assume throughout that $\Omega$ is a bounded domain (open
and connected) in $\R^n$ with Lipschitz boundary $\partial\Omega$.

\begin{definition}
  \label{def:3.1}
Given a measurable set $E\subset \R^n$ with positive measure and
given $u\in L^1(E)$ we denote its average
\begin{equation}
  \label{eq:3.1}
  u_E=\frac{1}{\mu_n(E)}\int_E u(x)d\mu_n(x) \ .
\end{equation}
\end{definition}
In the sequel we shall use this notation also for vector and
matrix valued functions, in particular we shall denote
\begin{equation}
  \label{eq:3.2}
  x_E=\frac{1}{\mu_n(E)}\int_E xd\mu_n(x) \ ,
\end{equation}
the center of mass of $E$.

We shall use analogous notation for averages on
$(n-1)$-dimensional Lipschitz surfaces $\Gamma\subset\R^n$
\begin{equation}
  \label{eq:3.3}
  u_\Gamma=\frac{1}{\mu_{n-1}(\Gamma)}\int_\Gamma
  u(x)d\mu_{n-1}(x) \ ,
\end{equation}
\begin{equation}
  \label{eq:3.4}
  x_\Gamma=\frac{1}{\mu_{n-1}(\Gamma)}\int_\Gamma
  xd\mu_{n-1}(x) \ .
\end{equation}
Here $\mu_{n-1}$ denotes the $(n-1)$-dimensional Lebesgue measure
on $\Gamma$. When no ambiguity occurs we shall also denote
$|E|=\mu_n(E)$, $|\Gamma|=\mu_{n-1}(\Gamma)$.

We recall the well-known Poincar\'e inequality.
\begin{theo}
  \label{theo:Poincare_classica}
There exists $Q>0$ such that
\begin{equation}
  \label{eq:3.Poincare_classica}
  \|u-u_\Omega\|_{L^2(\Omega)}\leq Q\|\nabla u\|_{L^2(\Omega)} \ ,
  \quad \hbox{for every }u\in W^{1,2}(\Omega) \ .
\end{equation}
\end{theo}
\begin{proof} See, for instance, \cite[Theorem 3.6.5]{l:m}. For a
quantitative evaluation of the constant $Q$ in terms of the
Lipschitz character of $\Omega$ we refer to \cite[Proposition
3.2]{l:amr2}
\end{proof}
We also recall the following generalized version, which is a
special case of a theorem due to Meyers \cite[Theorem 1]{l:mey}.
\begin{theo}
  \label{theo:Poincare_generica}
For every $\varphi\in(W^{1,2}(\Omega))^*$ such that
  $\varphi(1)=1$ we have
\begin{multline}
  \label{eq:3.Poincare_generica}
  \|u-\varphi(u)\|_{W^{1,2}(\Omega)}\leq
  \left(1+\|\varphi\|_{(W^{1,2}(\Omega))^*}|\Omega|^{\frac{1}{2}}\right)\sqrt{1+
Q^2}
  \|\nabla u\|_{L^2(\Omega)} \ ,\\
  \qquad \qquad \qquad \hbox{for every }u\in W^{1,2}(\Omega) \ .
\end{multline}
\end{theo}
\begin{rem}
  \label{rem:3.1}
Here we have chosen as the $W^{1,2}(\Omega)$-norm the expression
\begin{equation}
  \label{eq:3.norma}
  \|u\|_{W^{1,2}(\Omega)}=
  \left(\int_\Omega u^2+|\nabla u|^2\right)^{\frac{1}{2}} \ .
\end{equation}
Correspondingly, $\|\cdot\|_{(W^{1,2}(\Omega))^*}$ denotes the
dual norm.
\end{rem}
\begin{proof} [Proof of Theorem \ref{theo:Poincare_generica}]
This is in fact a straightforward consequence of Theorem
\ref{theo:Poincare_classica} and of Lemma \ref{lem:2.1}.
\end{proof}

\begin{ex}
\label{ex:Poincare_E}
Let $E$ be any measurable subset of $\Omega$
such that $|E|>0$. By choosing $\varphi(u)=u_E$, and computing
\begin{equation*}
  |\varphi(u)|\leq \frac{1}{\sqrt{|E|}}\|u\|_{L^{2}(\Omega)}\leq
  \frac{1}{\sqrt{|E|}}\|u\|_{W^{1,2}(\Omega)} \ , \quad \hbox{for
  every }u\in W^{1,2}(\Omega) \ ,
\end{equation*}
we obtain
\begin{equation}
\label{eq:3.Poincare_E}
  \|u-u_E\|_{W^{1,2}(\Omega)}\leq
  \left(1+\left(\frac{|\Omega|}{|E|}\right)^{\frac{1}{2}}\right)
  \sqrt{1+Q^2}
  \|\nabla u\|_{L^2(\Omega)} \ ,\quad \hbox{for
  every }u\in W^{1,2}(\Omega) \ .
\end{equation}
\end{ex}

\begin{ex}
\label{ex:Poincare_Gamma} Let $\Gamma$ be an open portion of
$\partial\Omega$. Let us denote by $\gamma(u)$ the trace on
$\Gamma$ of any $u\in W^{1,2}(\Omega)$, and let us choose
\begin{equation*}
  \varphi(u)=(\gamma(u))_\Gamma \ .
\end{equation*}
We evaluate
\begin{equation*}
  |\varphi(u)|\leq \frac{1}{\sqrt{|\Gamma|}}\|\gamma(u)\|_{L^{2}(\Gamma)}\leq
  \frac{C_\Gamma}{\sqrt{|\Gamma|}}\|u\|_{W^{1,2}(\Omega)} \ , \quad \hbox{for
  every }u\in W^{1,2}(\Omega) \ ,
\end{equation*}
where $C_\Gamma$ is the constant in the inequality for the trace
imbedding \mbox{$\gamma: W^{1,2}(\Omega)\rightarrow L^2(\Gamma)$}
and we obtain
\begin{multline}
\label{eq:3.Poincare_Gamma}
  \|u-(\gamma(u))_\Gamma\|_{W^{1,2}(\Omega)}\leq
  \left(1+C_\Gamma\left(\frac{|\Omega|}{|\Gamma|}\right)^{\frac{1}{2}}\right)
  \sqrt{1+Q^2}
  \|\nabla u\|_{L^2(\Omega)} \ ,\\
  \qquad \qquad\qquad\hbox{for
  every }u\in W^{1,2}(\Omega) \ .
\end{multline}
\end{ex}

\section{A higher order Poincar\'e inequality} \label{sec:
Poincare_H2}

We treat an analogue of Theorem \ref{theo:Poincare_generica}
suitable for functions in $W^{2,2}(\Omega)$.

We denote by $\mathcal{A}$ the space of affine functions on
$\Omega$
\begin{equation} \label{eq:4.affini}
  \mathcal{A}=\{u(x)=a+b\cdot x \ |\ a\in \R, b\in \R^n\} \ .
\end{equation}
The following Theorem is again a special case of Theorem 1 in
\cite{l:mey}. We report it here with a slightly different
expression of the constant on the right hand side.
\begin{theo}
  \label{theo:Poincare_generica_H2}
For any bounded operator $T:W^{2,2}(\Omega)\rightarrow\mathcal{A}$
such that $Tu=u$ for every $u\in \mathcal{A}$, we have
\begin{equation}
\label{eq:4.Poincare_generica_H2}
  \|u-Tu\|_{W^{2,2}(\Omega)}\leq
  (1+\|T\|)(1+Q^2+Q^4)^{\frac{1}{2}}
  \|\nabla^2 u\|_{L^2(\Omega)} \ ,\quad \hbox{for
  every }u\in W^{2,2}(\Omega) \ .
\end{equation}
Here $Q$ is the constant for the Poincar\'e inequality
\eqref{eq:3.Poincare_classica}.\end{theo}
\begin{proof}
The proof is immediate, by a repeated application of Poincar\'e
inequality \eqref{eq:3.Poincare_classica} and using Lemma
\ref{lem:2.1}.
\end{proof}
\begin{ex}
   \label{ex:Poincare_E_H2}
   In analogy to Example \ref{ex:Poincare_E}, given
   $E\subset\Omega$ measurable, with $|E|>0$, we pose
\begin{equation*}
  Tu=u_{E}+(\nabla u)_E\cdot(x-x_E) \ .
\end{equation*}
Note that $T$ is uniquely determined by the restriction of $u$ to
$E$. In fact if $u_1$, $u_2\in W^{1,2}(\Omega)$ are such that
$u_{1}|_{E}=u_{2}|_{E}$ then $\nabla(u_1-u_2)=0$ almost everywhere
in $E$ (see \cite[Lemma 7.7]{l:gt}).

Clearly $T$ is the identity on $\mathcal{A}$ and we compute
\begin{equation*}
  \|Tu\|_{W^{2,2}(\Omega)}^2\leq
  2\frac{|\Omega|}{|E|}(1+(\hbox{diam}\ \Omega)^2)\|u\|_{W^{1,2}(\Omega)}^2
\end{equation*}
and consequently we obtain
\begin{multline}
\label{eq:4.Poincare_E_H2}
  \|u-(u_{E}+(\nabla u)_E\cdot(x-x_E))\|_{W^{2,2}(\Omega)}\leq\\
  \leq\left(1+\left(2
  \frac{|\Omega|}{|E|}(1+(\hbox{diam}\ \Omega)^2)\right)^{\frac{1}{2}}\right)
  (1+Q^2+Q^4)^{\frac{1}{2}}\|\nabla^2 u\|_{L^2(\Omega)} \ ,\\
  \qquad\qquad\qquad\hbox{for every } u\in W^{2,2}(\Omega) \ .
\end{multline}
\end{ex}

\begin{ex}
   \label{ex:Poincare_balls_H2}
   In certain instances, it is useful to have an estimate  of the
   type \eqref{eq:4.Poincare_E_H2} when $\Omega$ and $E$ are
   concentric balls, say $\Omega=B_1(0)$, $E=B_\rho(0)$,
   $0<\rho<1$. In this case we evaluate
\begin{equation*}
  \|T\|\leq
  \sqrt{\frac{n+3}{n+2}}\left(\frac{1}{\rho}\right)^{\frac{n}{2}}
  \ .
\end{equation*}
Consequently
\begin{multline}
\label{eq:Poincare_balls_H2}
  \|u-Tu\|_{W^{2,2}(B_1(0))}\leq
  \left(1+\sqrt{\frac{n+3}{n+2}}\left(\frac{1}{\rho}\right)^{\frac{n}{2}}
  \right)(1+Q^2+Q^4)^{\frac{1}{2}}
  \|\nabla^2 u\|_{L^2(B_1(0))} \ ,\\
  \qquad \qquad \hbox{for
  every }u\in W^{2,2}(B_1(0)) \ .
\end{multline}
This specific estimate turns out to be useful in an inverse
problem for elastic plates, \cite{l:mrv}.
\end{ex}
\begin{ex}
   \label{ex:Poincare_Gamma_H2}
   We note that in general it is not possible to mimic Example
   \ref{ex:Poincare_Gamma} for a Poincar\'e type inequality in
   $W^{2,2}(\Omega)$ when some functionals of the trace of $u$ \emph{alone} on
   an open portion $\Gamma$ of the boundary are known.

   For example, if one considers
\begin{equation*}
 \Omega=\{x\in \R^n\ |\ 0<x_i<1, i=1,\ldots,n\}
\end{equation*}
and
\begin{equation*}
 \Gamma=\{x\in \R^n\ |\ x_n=0, 0<x_i<1, i=1,\ldots,n-1\} \ ,
\end{equation*}
then the function $u=x_n$ has zero trace on $\Gamma$ and
$\|\nabla^2 u\|_{L^2(\Omega)}=0$. Hence its $W^{2,2}(\Omega)$-norm
is not dominated by $\|\nabla^2 u\|_{L^2(\Omega)}$ and any
functional of its trace on $\Gamma$.

In fact some additional assumptions, or data, are necessary.

If we admit that $\Gamma$ may be flat, additional pieces of
information are needed, see \cite{l:mey}. For instance, if also
the trace of $\nabla u$ on $\Gamma$ is available, then one can
consider
\begin{equation*}
Tu=(\gamma(u))_\Gamma+(\gamma(\nabla u))_\Gamma\cdot(x-x_\Gamma)
\end{equation*}
and obtain
\begin{multline}
\label{eq:4.Poincare_H2_no_flat}
\|u-(\gamma(u))_\Gamma-(\gamma(\nabla
u))_\Gamma\cdot(x-x_\Gamma)\|_{W^{2,2}(\Omega)}\leq
(1+\|T\|)(1+Q^2+Q^4)^{\frac{1}{2}} \|\nabla^2 u\|_{L^2(\Omega)} \ ,\\
  \qquad\qquad\qquad\hbox{for every }u\in W^{2,2}(\Omega) \ .
\end{multline}
And we compute
\begin{equation}
  \label{eq:4.Poincare_H2_no_flat_norm}
\|T\|\leq C_\Gamma\left(2\frac{|\Omega|}{|\Gamma|}
(1+(\hbox{diam}\ \Omega)^2)\right)^{\frac{1}{2}}
\end{equation}
where $C_\Gamma$ is the constant in the inequality for the trace
imbedding \mbox{$\gamma: W^{1,2}(\Omega)\rightarrow L^2(\Gamma)$}.

Let us assume, instead, that $\Gamma\subset\partial\Omega$ is not
flat, that is, it is not a portion of an hyperplane. In this case
we obtain Theorem \ref{theo:funzionale_traccia} below. In order to
state it, we first need some preparation.

 We choose a
reference system whose origin lies on the center of mass of
$\Gamma$, that is, we assume
\begin{equation}
\label{eq:4.baricentro_origine}
  x_\Gamma=0 \ .
\end{equation}
\end{ex}

\begin{lem}\label{traccia di A}
Assume that $\Gamma$ is not a portion of an hyperplane. The
restriction to $\mathcal{A}$ of the trace imbedding
$\gamma:W^{1,2}(\Omega)\rightarrow L^2(\Gamma)$ is one to one.
\end{lem}
\begin{proof}First of all we observe that, in
view of \eqref{eq:4.baricentro_origine},  $\Gamma$ is not a
portion of an hyperplane if and only if there exists $n$ linearly
independent vectors $v_1,\ldots,v_n\in\Gamma$. Thus, if $l(x)=a +b
\cdot x \in \mathcal{A}$ and $\gamma(l) = 0$, then $a = (
\gamma(l))_{\Gamma} = 0$ and $b \in \mathbb{R}^n$ must satisfy
\begin{equation*}
b \cdot v_i = \gamma(l)(v_i) - a = 0 \ , \ \text{ for every } i=1,
\ldots , n \ .
\end{equation*}
Consequently $b=0$ and hence, $l=0$.
\end{proof}

\begin{rem}
Observe that, being $\Omega$ bounded, any sufficiently large
portion $\Gamma$ of the boundary shall not be flat.\end{rem}

Let us denote by $e:\gamma(\mathcal{A})\rightarrow \mathcal{A}$
the inverse of $\gamma|_{\mathcal{A}}$, and by
\mbox{$\pi:L^2(\Gamma) \rightarrow \gamma(\mathcal{A})$} the
orthogonal projection onto $\gamma(\mathcal{A})$ with respect to
the $L^2(\Gamma)$ inner product. Note that an explicit expression
of $\pi$ is easily obtained. In fact, for every $v\in
L^2(\Gamma)$, $\pi(v)$ is determined as $\pi(v)=\gamma(l)$, where
$l\in\mathcal{A}$ is the minimizer of the finite dimensional least
squares problem
\begin{equation*}
\min \left\{ \int_{\Gamma}(v - l)^2  | \ l \in \mathcal{A}
\right\} \ .
\end{equation*}
We also introduce $\tau:L^2(\Gamma)\rightarrow \mathcal{A}$ as the
composition $\tau = e\circ \pi$.
\begin{theo}
  \label{theo:funzionale_traccia}
  Assume that $\Gamma$ is not a portion of an hyperplane. We have
\begin{multline}
\label{eq:4.funzionale_traccia}
  \|u-\tau(\gamma(u))\|_{W^{2,2}(\Omega)}\leq(1+C_\Gamma\|e\|)
  (1+Q^2+Q^4)^{\frac{1}{2}}\|\nabla^2 u\|_{L^2(\Omega)} \ ,\\
  \qquad\qquad\qquad\hbox{for every }u\in W^{2,2}(\Omega) \ .
\end{multline}
Here $\|e\|$ denotes the ${\cal
L}(L^2(\Gamma),W^{2,2}(\Omega))$-norm of $e$, $C_\Gamma$ is the
constant in the inequality for the trace imbedding \mbox{$\gamma:
W^{1,2}(\Omega)\rightarrow L^2(\Gamma)$} and $Q$ is the constant
for the Poincar\'e inequality \eqref{eq:3.Poincare_classica}.
\end{theo}
\begin{rem}\label{rem:e}
Observe that the calculation of the norm $\|e\|$ reduces to
solving a finite dimensional eigenvalue problem, in fact
\begin{equation*}
\|e\|^2 = \max \left\{ \ \frac{\int_{\Omega}(a + b \cdot x)^2 +
|b|^2}{\int_{\Gamma}(a + b \cdot x)^2 }\ | \ a\in \mathbb{R}\ , \
b \in \mathbb{R}^n  \ , \ a^2+|b|^2 > 0 \right\} \ .
\end{equation*}
\end{rem}
As a consequence of Theorem \ref{theo:funzionale_traccia} we
obtain the following Corollary.

\begin{cor}
   \label{rem:corollario_funzionale_traccia}
  If  $\Gamma$ is not a portion of an hyperplane, then we have
\begin{equation}
\label{eq:4.corollario_funzionale_traccia}
\|u\|_{W^{2,2}(\Omega)}\leq C_1\|\nabla^2 u\|_{L^2(\Omega)}+C_2
\|\gamma(u)\|_{L^2(\Gamma)} \ ,
  \quad\hbox{for every }u\in W^{2,2}(\Omega) \ ,
\end{equation}
where $C_1$, $C_2>0$ are constants depending on $\Omega$ and
$\Gamma$ only.
\end{cor}

\begin{proof} [Proof of Theorem
\ref{theo:funzionale_traccia}] It suffices to apply Theorem
\ref{theo:Poincare_generica_H2} with $T=\tau \circ \gamma$ and to
observe that, since
$\|\pi\|_{\mathcal{L}(L^2(\Gamma),L^2(\Gamma))}=1$, then we have $
\|T\| \leq C_{\Gamma}\|\tau\|\leq C_{\Gamma}\|e\|$.

\end{proof}

\section{Korn type inequalities} \label{sec:
Korn} We recall the Korn inequality (of second kind). Given $u\in
W^{1,2}(\Omega;\R^n)$ we denote
\begin{equation*}
\widehat \nabla u=\frac{1}{2}(\nabla u+(\nabla u)^T) \ .
\end{equation*}
\begin{theo}
  \label{theo:korn_classica}
  There exists $K>0$ such that
\begin{equation}
  \label{eq:5.korn_classica}
\left\|\nabla u-\frac{1}{2}(\nabla u-(\nabla
u)^T)_\Omega\right\|_{L^2(\Omega)}\leq K\|\widehat \nabla
u\|_{L^2(\Omega)} \ , \quad\hbox{for every }u\in
W^{1,2}(\Omega;\R^n) \ .
\end{equation}
\end{theo}
\begin{proof}
See for instance \cite{l:f}.
\end{proof}
We denote by $\hbox{Skew}^n=\{M\in \M^{n \times n}\ |\ M^T=-M\}$
the class of skew-symmetric $n \times n$ matrices and by ${\cal
R}=\{b+Mx\ |\ b\in \R^n, M\in \hbox{Skew}^n\}$ the linear space of
infinitesimal rigid displacements.
\begin{theo}
  \label{theo:korn_generica}
For any bounded linear operator $T:W^{1,2}(\Omega;\R^n)\rightarrow
{\cal R}$ such that $Tu=u$ for every $u\in {\cal R}$, we have
\begin{equation}
  \label{eq:5.korn_generica}
\|u-Tu\|_{W^{1,2}(\Omega)}\leq
(1+\|T\|)K(1+Q^2)^{\frac{1}{2}}\|\widehat \nabla u\|_{L^2(\Omega)}
\ , \quad\hbox{for every }u\in W^{1,2}(\Omega;\R^n) \ .
\end{equation}
Here $K$ is the constant appearing in \eqref{eq:5.korn_classica}
and $Q$ is the constant in the Poincar\'e inequality
\eqref{eq:3.Poincare_classica}.
\end{theo}

\begin{proof}
Set $P:W^{1,2}(\Omega;\R^n)\rightarrow{\cal R}$ as
\begin{equation*}
Pu=u_\Omega+\frac{1}{2}(\nabla u-(\nabla u)^T)_\Omega (x-x_\Omega)
\ .
\end{equation*}
By \eqref{eq:3.Poincare_classica} and \eqref{eq:5.korn_classica}
we obtain
\begin{equation*}
\|u-Pu\|_{W^{1,2}(\Omega)}\leq K(1+Q^2)^{\frac{1}{2}}\|\widehat
\nabla u\|_{L^2(\Omega)} \ .
\end{equation*}
The thesis follows {}from  Lemma \ref{lem:2.1} by choosing
\begin{equation*}
X=W^{1,2}(\Omega;\R^n) \ ,
\end{equation*}
\begin{equation*}
Y=L^{2}(\Omega;\M^{n \times n}) \ ,
\end{equation*}
\begin{equation*}
Lu=\widehat \nabla u \ ,
\end{equation*}
\begin{equation*}
X_0={\cal R} \ .
\end{equation*}
\end{proof}
\begin{ex}
\label{ex:korn_E} Given $E\subset\Omega$ measurable with $|E|>0$
and setting $Tu=u_E+\frac{1}{2}(\nabla u-(\nabla u)^T)_E (x-x_E)$,
we compute
\begin{equation*}
  \|Tu\|_{W^{1,2}(\Omega)}^2\leq
  2\frac{|\Omega|}{|E|}(1+(\hbox{diam}\ \Omega)^2)\|u\|_{W^{1,2}(\Omega)}^2 \ .
\end{equation*}
Consequently, we obtain
\begin{multline}
\label{eq:5.korn_E}
  \left\|u-u_E-\frac{1}{2}(\nabla u-(\nabla
  u)^T)_E (x-x_E)\right\|_{W^{1,2}(\Omega)}\leq\\
  \leq\left(1+\left(2\frac{|\Omega|}{|E|}
  (1+(\hbox{diam}\ \Omega)^2)\right)^{\frac{1}{2}}\right)K(1+Q^2)^{\frac{1}{2}}
  \|\widehat\nabla u\|_{L^2(\Omega)} \ ,\\
  \qquad\qquad\qquad \hbox{for
  every }u\in W^{1,2}(\Omega;\R^n) \ .
\end{multline}
\end{ex}

\begin{ex}
\label{ex:korn_balls}
   When $\Omega$ and $E$ are
   concentric balls, say $\Omega=B_1(0)$, $E=B_\rho(0)$,
   $0<\rho<1$. In this case we evaluate
\begin{equation*}
  \|T\|\leq
  \sqrt{\frac{n+3}{n+2}}\left(\frac{1}{\rho}\right)^{\frac{n}{2}}.
\end{equation*}
Consequently
\begin{multline}
\label{eq:Poincare_balls_H2}
  \|u-Tu\|_{W^{1,2}(B_1(0))}\leq
  \left(1+\sqrt{\frac{n+3}{n+2}}\left(\frac{1}{\rho}\right)^{\frac{n}{2}}
  \right)K(1+Q^2)^{\frac{1}{2}}
  \|\widehat{\nabla} u\|_{L^2(B_1(0))} \ ,\\
  \qquad \qquad \hbox{for
  every }u\in W^{1,2}(B_1(0);\R^n) \ .
\end{multline}
This bound improves the one in \cite[Lemma 3.5]{l:amr1} in two
respects. It applies to any dimension and the exponent of
$\frac{1}{\rho}$ is diminished.
\end{ex}

\begin{ex}
\label{ex:korn_Gamma} We examine a version of Korn inequality when
constraints on $u$ are taken on its trace on an open portion
$\Gamma$ of $\partial\Omega$. The following Lemma will be useful.
\end{ex}

\begin{lem}\label{traccia di R}
 The
restriction to $\mathcal{R}$ of the trace imbedding
\mbox{$\gamma:W^{1,2}(\Omega;\R^n)\rightarrow L^2(\Gamma;\R^n)$}
is one to one.
\end{lem}
\begin{proof}
With no loss of generality we assume $x_{\Gamma}=0$. Being
$\Gamma$ an $(n-1) - $ dimensional hypersurface, there exist
$v_1,\ldots,v_{n-1}\in\Gamma$ which are linearly independent. We
show that, given $r\in \mathcal{R}$ such that $\gamma(r)=0$ we
have $r=0$. In fact we have $r(0)=\gamma(r)_{\Gamma}=0$ and also
$r(v_i)=\gamma(r)(v_i)=0$ for every $i=1, \ldots ,n-1$. It is well
known that such constraints imply $r=0$. We provide a proof for
the sake of completeness. Let us write $r(x)= b + Mx$ with $b \in
\R^n$ and $M \in \hbox{Skew}^n$. We have $b=r(0)=0$, hence it
remains to prove that the conditions $r(v_i)=0 \ , i=1, \ldots
,n-1$ imply $M=0$. Up to a rotation in the reference system, we
can assume that the $n-$th component of the vectors $v_i$,
$i=1,\ldots,n-1$, is zero. Thus we have
\begin{equation*}
\left(v_1|\cdots|v_{n-1}|0 \right)= \left(
 \begin{array}{ccc | c}
   &&&0\\
   & W && \vdots\\
   &&& 0 \\
   \hline
   0&\cdots&0&0
  \end{array}
\right) \ ,
\end{equation*}
where $W$ is a nonsingular $(n-1)\times(n-1)$ matrix. We have
\begin{equation*}
Mv_i = r(v_i) = 0 \ , \text{ for every } i=1,\ldots,n-1 \ .
\end{equation*}
Consequently
\begin{equation*}
M \left(
 \begin{array}{ccc | c}
   &&&0\\
   & W && \vdots\\
   &&& 0 \\
   \hline
   0&\cdots&0&0
  \end{array}
\right)= 0\ ,
\end{equation*}
and also
\begin{equation*}
M \left(
 \begin{array}{ccc | c}
   &&&0\\
   & I && \vdots\\
   &&& 0 \\
   \hline
   0&\cdots&0&0
  \end{array}
\right) = M \left(
 \begin{array}{ccc | c}
   &&&0\\
   & W && \vdots\\
   &&& 0 \\
   \hline
   0&\cdots&0&0
  \end{array} \right) \left(
  \begin{array}{ccc | c}
   &&&0\\
   & W^{-1} && \vdots\\
   &&& 0 \\
   \hline
   0&\cdots&0&0
  \end{array}
\right) = 0\ ,
\end{equation*}
where $I$ is the $(n-1)\times(n-1)$ identity matrix. It follows
that the entries below the diagonal of $M$ are all zero and,
consequently, $M=0$.
\end{proof}
Let us denote by $E:\gamma(\mathcal{R})\rightarrow \mathcal{R}$
the inverse of $\gamma|_{\mathcal{R}}$, and by
\mbox{$\Pi:L^2(\Gamma;\R^n) \rightarrow \gamma(\mathcal{R})$} the
orthogonal projection onto $\gamma(\mathcal{R})$ with respect to
the $L^2(\Gamma;\R^n)$ inner product. We also introduce
$\rho:L^2(\Gamma;\R^n)\rightarrow {\cal R}$ as the composition
$\rho = E\circ \Pi$.
\begin{theo}
  \label{theo:funzionale_traccia_korn}
   Given an open portion $\Gamma$ of $\partial\Omega$, we have

\begin{multline}
\label{eq:5.funzionale_traccia_korn}
  \|u-\rho(\gamma(u))\|_{W^{1,2}(\Omega)}\leq(1+C_\Gamma\|E\|)
  K(1+Q^2)^{\frac{1}{2}}\|\widehat\nabla u\|_{L^2(\Omega)} \ ,\\
  \qquad\qquad\qquad\hbox{for every }u\in W^{1,2}(\Omega;\R^n) \  .
\end{multline}
Here $\|E\|$ denotes the ${\cal
L}(L^2(\Gamma;\R^n),W^{1,2}(\Omega;\R^n))-$norm of $E$,
$C_\Gamma$ is the constant in the inequality for the trace
imbedding \mbox{$\gamma: W^{1,2}(\Omega)\rightarrow L^2(\Gamma)$}
and $Q$ is the constant for the Poincar\'e inequality
\eqref{eq:3.Poincare_classica}.
\end{theo}
\begin{rem}
Arguing as in Remark \ref{rem:e}, also the calculation of the norm $\|E\|$ reduces to solving a finite dimensional eigenvalue problem.
\end{rem}

The proof is straightforward, as well as the one of the following
Corollary.
\begin{cor}
   \label{cor:corollario_funzionale_traccia_korn}
   There exist $C_1$, $C_2>0$, depending only on
   $\Omega$ and $\Gamma$, such that
\begin{equation}
\label{eq:5.corollario_funzionale_traccia_korn}
\|u\|_{W^{1,2}(\Omega)}\leq C_1\|\widehat\nabla
u\|_{L^2(\Omega)}+C_2 \|\gamma(u)\|_{L^2(\Gamma)} \ ,
  \quad\hbox{for every }u\in W^{1,2}(\Omega;\R^n) \ .
\end{equation}
\end{cor}

\noindent \textbf{Acknowledgement.} The authors wish to thank an anonymous referee for providing fundamental bibliographical information which
was missing in a previous version of this paper.


\begin{thebibliography}{a}

\bibitem{l:amr1}
G.~Alessandrini, A.~Morassi and E.~Rosset.
\newblock {Detecting an inclusion in an elastic body by boundary
measurements}.
\newblock {\em SIAM J. Math. Anal.} \textbf{33} (6) (2002) 1247--1268.

\bibitem{l:amr2}
G.~Alessandrini, A.~Morassi and E.~Rosset.
\newblock {Detecting cavities by electrostatic boundary measurements}.
\newblock {\em Inverse Problems} \textbf{18} (2002) 1333--1353.

\bibitem{l:f}
K.O.~Friedrichs.
\newblock {On the boundary value problems of the theory of elasticity and
Korn's inequality}.
\newblock {\em Annals of Math.} \textbf{48} (1947) 441--471.

\bibitem{l:gt}
D.~Gilbarg and N.S.~Trudinger.
\newblock {\em Elliptic partial
differential equations of second order}.
\newblock Springer, New York, 1983.

\bibitem{l:ko}
V.A.~Kondrat'ev and O.A.~Oleinik.
\newblock {On the dependence of the constant in Korn's inequality
on parameters characterizing the geometry of the region}.
\newblock {\em Russian Math. Surveys} \textbf{44} (1989) 187--195.

\bibitem{l:mey}
N.G.~Meyers.
\newblock {Integral inequalities of Poincar\'{e} and Wirtinger type}.
\newblock {\em Arch. Rational Mech. Anal.} \textbf{68} (2) (1978)
113--120.

\bibitem{l:meyz}
N.G.~Meyers and W.P.~Ziemer.
\newblock {Integral inequalities of Poincar\'{e} and Wirtinger type for BV functions}.
\newblock {\em Amer. J. Math.} \textbf{99} (6) (1977)
1345--1360.

\bibitem{l:mr}
A.~Morassi and E.~Rosset.
\newblock {Uniqueness and stability
in determining a rigid inclusion in an elastic body}.
\newblock Submitted for publication.

\bibitem{l:mrv}
A.~Morassi, E.~Rosset and S.~Vessella.
\newblock {Size estimates for inclusions in an elastic plate by boundary measurements}.
\newblock Submitted for publication.

\bibitem{l:m}
C.B.~Morrey.
\newblock {\em Multiple Integrals in the Calculus of Variations}.
\newblock {Springer-Verlag}, New York, 1966.

\bibitem{l:z}
W.P.~Ziemer.
\newblock {\em Weakly Differentiable Functions}.
\newblock {Springer-Verlag}, New York, 1989.


\end{thebibliography}
\end{document}